\documentclass[]{article}
\usepackage{amsthm,amsmath,amssymb}
\usepackage[OT1]{fontenc}
\usepackage[colorlinks,citecolor=blue,urlcolor=blue]{hyperref}
\usepackage{mathtools, bm, bbold}
\usepackage{parskip}
\usepackage{enumitem}
\usepackage[ruled]{algorithm2e}
\usepackage{authblk}
\usepackage{tabu,booktabs,float,tabularx}
\usepackage{subcaption,caption}
\usepackage{siunitx}

\usepackage[noabbrev,capitalize]{cleveref}
\usepackage[maxnames=999,giveninits=true]{biblatex}
\usepackage{xcolor}

\bibliography{assets/all}
\DeclareFieldFormat{titlecase}{\MakeSentenceCase*{#1}}

\newtheorem{theorem}{Theorem}[section]

\newtheorem{dataset}[theorem]{Setting}

\renewcommand{\phi}{\ensuremath{\varphi}}
\renewcommand{\epsilon}{\ensuremath{\varepsilon}}
\title{Local moment matching with Erlang mixtures under automatic roughness penalization}
\author[1]{Oskar Laverny \thanks{E-mail address: \href{mailto:oskar.laverny@univ-amu.fr}{oskar.laverny@univ-amu.fr}}}
\author[2,3]{Philippe Lambert}
\affil[1]{Aix Marseille Univ, Inserm, IRD, SESSTIM, Sciences Economiques \& Sociales de la Santé \& Traitement de l’Information Médicale, ISSPAM, Marseille, France.}
\affil[2]{Institut de Statistique, Biostatistique et Sciences Actuarielles (ISBA), Université catholique de Louvain, Voie du Roman Pays 20, B-1348 Louvain-la-Neuve, Belgium}
\affil[3]{Institut de Mathématique, Université de Liège, Allée de la Découverte 12, B-4000 Liège, Belgium.} 

\begin{document}
  \maketitle
  \begin{abstract}We consider the class of Erlang mixtures for the task of density estimation on the positive real line when the only available information is given as local moments, a histogram with potentially higher order moments in some bins. By construction, the obtained moment problem is ill-posed and requires regularization. Several penalties can be used for such a task, such as a lasso penalty for sparsity of the representation, but we focus here on a simplified roughness penalty from the P-splines literature. We show that the corresponding hyperparameter can be selected without cross-validation through the computation of the so-called effective dimension of the estimator, which makes the estimator practical and adapted to these summarized information settings. The flexibility of the local moments representations allows interesting additions such as the enforcement of Value-at-Risk and Tail Value-at-Risk constraints on the resulting estimator, making the procedure suitable for the estimation of heavy-tailed densities.
\end{abstract}


\clearpage 
\section{Introduction}

Let $X_1,...,X_n$ be a $n$-sample of a positive real random variable $X$, and let $\mathcal B = \left(B_1,...,B_J\right)$ be a finite partition of $\mathbb R_+$, where $B_j = [b_{j-1}, b_j[$ are called \emph{bins}.
Let furthermore $\bm k \in \mathbb N^J$ be a vector of integers. We consider that the only information we have about the sample is the set of empirical moments given by: 

\begin{equation}\label{eqn:pihatandmuhat}
\begin{aligned}
  \hat{\bm \pi} &= \left\{\hat{\pi}_{j} = \frac{1}{N}\sum_{i=1}^N \mathbb 1_{X_i \in B_j}: j \in 1,...,J\right\}\text{ and}\\
  \hat{\bm \mu} &= \left\{\hat{\mu}_{j,k} = \frac{1}{N} \sum_{i=1}^N X_i^k\mathbb 1_{X_i \in B_j} : j \in 1,...J, k \in 1,...,k_j\right\}.
\end{aligned}
\end{equation}

The zeroth moments are denoted separately strictly for ease of notations. In an insurance or reinsurance context, such binned information may arise from attritionnal versus large cutoff in a loss random variable, and the lack of detailed information (only a few moments are available) may typically result from confidentiality concerns. 

Indeed, the General Data Protection Regulation (GDPR) has raised important confidentiality concerns in the insurance industry. The regulation imposes strict rules on the collection, use, and storage of personal data, including data related to health, which is often used in underwriting and claims management. The GDPR's requirements for confidentiality in insurance pose significant challenges for insurers, who must balance the need to collect and use personal data with the need to protect individuals' privacy. 
Several studies have examined the impact of this regulation. While providing an implementation guide, \cite{liapakis2018gdpr} notes that the process could be used by the industry as an opportunity for fine-tuning the data management. On the other hand, \cite{dexe2020empirical} investigates the way insurer respond to the the right to explanation and the need for transparency in premium computations.

Because of these concerns, it is not uncommon for losses data to be transmitted in a summarized state alike \cref{eqn:pihatandmuhat}, with rather small integers $J,k_1,...k_J$, to increase the protection of personal data by not transmitting the full losses and/or premium information.
A representative example of such a problem drawn from a \texttt{LogNormal}$(0,0.5)$\footnote{We denote all probability distributions by $\texttt{Name}(\text{Parameters})$. Their definitions are all detailed in \cref{apx:keydists}.} is given in \cref{tab:exemple_1}.

\begin{table}[H]
	\caption{\label{tab:exemple_1}Example of a local moment problem drawn from a \texttt{LogNormal}$(0,0.5)$ sample. In this example, the three first lines gives estimated boxed moments. The fourth line is a little more involved: it prescribes a Value-At-Risk $b_{J-1}$ and a Tail-Value-at-Risk $\hat{\mu}_{J-1,1}$ at the quantile level $1-\hat{\pi}_4$. This data structure is classical when modeling insurance losses, which are usually divided into \emph{attritional} and \emph{large} losses.}
	\centering
  \begin{tabularx}{\linewidth}{c|llllXXXX}
\toprule 
$j$ & $[b_{j-1},b_j[$ & $N_j$ & $k_j$ & $\pi_j$ & $\hat{\mu}_{j,1}$ & $\hat{\mu}_{j,2}$ & $\hat{\mu}_{j,3}$ & $\hat{\mu}_{j,4}$ \\
\midrule 
1 & $[0.000,0.948[$ & 375 & 4 & 0.500 & 0.332 & 0.235 & 0.175 & 0.136 \\
2 & $[0.948,1.885[$ & 300 & 4 & 0.400 & 0.526 & 0.719 & 1.017 & 1.488 \\
3 & $[1.885,3.332[$ & 67 & 4 & 0.089 & 0.206 & 0.485 & 1.167 & 2.874 \\
4 & $[3.332,\infty[$ & 8 & 1 & 0.011 & 0.048 &  &  &  \\
\bottomrule 
\end{tabularx}

\end{table}

The unavailability of higher order moments on the last line of \cref{tab:exemple_1} represents the fact that these moments might not be finite and are therefore usually not available. This flexibility of the representation allows handling heavy-tailed random variables through VaR/TVaR insights (or even only VaR). In what follows, we only assume that the zeroth moments in each box are entirely part of the available information (that is $k_j \ge 0$ for all $j \in J$).

From this restricted information, our goal is to estimate the density for the positive (and potentially heavy-tailed) random variable $X$. If the previous literature on such binned moment problems (see, e.g.,~\cite{lambertNonparametricDensityEstimation2023} or the related \cite{brockett1997actuarial}) dealt mainly with finitely supported random variables, we are here particularly interested in the estimation of the tail of distribution. Our proposition is to tackle this problem using the class of Erlang mixtures.

Erlang Mixtures have been widely used these past 30 years in various fields, including telecommunications, manufacturing, reliability analysis and insurance, in which they are used to model the distribution of losses, claim amounts, and other risk-related variables. Their flexibility allows them to capture a wide range of distributions, including skewed and heavy-tailed distributions that are commonly found in insurance data. This has led to their adoption in various areas of insurance, such as pricing, reserving, and risk management. See, e.g., \cite{leeModelingEvaluatingInsurance2010,leeModelingDependentRisks2012} for pioneering work on their estimation through Expectation-Maximization (EM) algorithms, in univariate or multivariate settings respectively. The multivariate case is commonly used to model the joint distribution of multiple risk-related variables, including the dependence structure, such as claim amounts, frequency, and severity among several lines of business. See \cite{gongErlangbasedMethodsModeling2014} for applications to ruin problems (and references therein), or \cite{porth2014credibility} for application to crop insurance pricing. More recently, \cite{guiFittingErlangMixture2018,guiFittingMultivariateErlang2021a} uses multivariate Erlang mixtures to fit censored and truncated multivariate data with a roughness penalty, again through a (modified) EM algorithm. This literature shows that Erlang mixtures are able to capture complex tail behaviors and that the model can provide valuable insights on the underlying risk factors.

Overall, univariate and multivariate Erlang mixtures offer a flexible and powerful framework for analyzing insurance data, which makes them a valuable tool for insurance practitioners and researchers alike.
The denseness results from \cite{tijms1994stochastic} and subsequent multivariate extensions are probably at the root of the popularity of this model. We are particularly interested in the work in~\cite{cossette2016moment}, which solve standard moment problems on the real line with Erlang mixtures, again with risk management applications, giving an alternative to the classical estimator from~\cite{johnson1989matching}. More precisely, \cite{johnson1989matching} and \cite{cossetteDependentRiskModels2018} solve our problem for a unique, fully covering bin (that is, $J=1$ and $B_1 = \mathbb R_+$), corresponding to a standard moment problem on the positive real line.

If our estimation problem is a little more involved than a standard moment problem, it is easy to see that such an inverse problem will require regularization. Therefore, we propose to leverage a penalization of the roughness of the density, as was already discussed in~\cite{guiFittingMultivariateErlang2021a}. However, two problems arise: 
\begin{itemize}
  \item The direct penalization of the continuous roughness of the density proposed in~\cite{guiFittingMultivariateErlang2021a} comes with a large computational cost. 
  \item In our low information settings, we do not have a full dataset and the calibration of the penalty hyperparameter is not possible through cross-validation (like in classical EM schemes) since holding out testing data would not be acceptable.
\end{itemize}

On the other hand, the penalization of roughness of a density is largely dealt with in the literature around approximation by splines, notably yielding the concept of P-splines. 

P-splines, or penalized B-splines, are a popular and flexible tool for fitting smooth curves and surfaces to data. They were introduced by Eilers and Marx~\cite{eilersFlexibleSmoothingBsplines1996a} as a way to overcome some limitations of traditional B-spline methods, such as the difficulty of selecting knots and the optimal degree of smoothing. P-splines use a penalty on the $r^{\text{th}}$-order differences of the B-spline coefficients to achieve smoothness, and this penalty can be easily controlled through a tuning parameter. See~\cite{eilers2021practical} for extensive literature review on P-splines. 

Through a Bayesian perspective, we leverage the P-spline literature to provide in this paper a solution to the two problems encountered with the roughness penalization of Erlang mixtures. The discussion will be organized as follows: \cref{sec:proposed_est_scheme} derives and motivates the proposed approximation scheme, \cref{sec:numerical_investigation} gives a few numerical illustrations of the method, including convergence study on simulated datasets, and \cref{sec:conclusion} concludes. 

\section{Proposed estimation scheme}\label{sec:proposed_est_scheme}

Recall that we only observe $\hat{\bm \pi}$ and $\hat{\bm \mu}$ from \cref{eqn:pihatandmuhat}. We also consider that we know $N$, the number of observations that were used to compute the histogram $\hat{\bm \pi}$ and the boxed moments $\hat{\bm \mu}$. Let's assume without loss of generality that $X$ has theoretical counterpart moments: 

\begin{equation}\label{eqn:piandmu}
\begin{aligned}
  \bm \pi &= \left\{\pi_{j} = \mathbb P\left(X \in B_j\right): j \in 1,...,J\right\}\text{ and}\\
  \bm \mu &= \left\{\mu_{j,k} = \mathbb E\left(X^k\mathbb 1_{X \in B_j}\right) : j \in 1,...J, k \in 1,...,k_j\right\}.
\end{aligned}
\end{equation}

We consider the underlying sample $X_1,...X_n$ to be independently and identically distributed. Under this assumption, we have that $N\hat{\bm \pi} \sim \texttt{Multinomial}(\bm\pi,N)$ (see \cref{apx:keydists} for distribution definitions). We then approximate, conditionally on $\hat{\bm \pi}$, the distribution of the vector of moments $\hat{\bm \mu}$ through the law of large numbers as $\hat{\bm \mu} \sim \texttt{Normal}(\bm\mu, \bm \Sigma/N)$, where the variance-covariance matrix $\bm \Sigma$ is defined by: $$\Sigma_{(j,k), (i,m)} = \mu_{j,k+m}\mathbb 1_{j = i} - \mu_{j,k}\mu_{i,m}.$$

We generally denote by $\ell$ log-density, neglecting additive constants, of random objects. In particular, the loglikelihood of our model can be expressed (neglecting additive constants) as
\begin{equation}\label{eqn:llh}
  \ell(\bm\pi,\bm\mu,\bm\Sigma) = \left\{\hat{\bm\pi}'\log(\bm\pi)\right\} + \left\{- \frac{1}{2} \log|\bm\Sigma| - \frac{1}{2} \lVert\bm \mu - \hat{\bm \mu}\rVert_{\Sigma}^2\right\},
\end{equation}
where brackets delimit the respective contribution of the multinomial distribution and the conditional multivariate normal. Here, $\log(\bm\pi)$ is intended componentwise, $\lvert \cdot\rvert$ denotes the determinant of a matrix and $\lVert \bm x \rVert_{\Sigma}^2 = \bm x' \bm \Sigma^{-1} \bm x$ denotes the squared Mahalanobis distance. It is clear from this expression that, under this Gaussian approximation for $\hat{\bm\mu}$, the triplet $(\bm\pi,\bm\mu,\bm\Sigma)$ is sufficient to evaluate the loglikelihood of our data, whatever the underlying parametric or non-parametric model we choose.

Due to the positive and potentially heavy-tailed nature of the random variable, we propose to model it through the class of Erlang Mixtures. 
Recall that the $\texttt{MixedGamma}(\nu,\theta)$ distribution, where $\theta \in \mathbb R_+$ is the scale and $\nu \in \mathcal P(\mathbb R_+)$ is the prior distribution of shapes, has density 
$$f_{\nu,s}(x) = \int_{\mathbb R_+} \frac{x^{\alpha-1}}{\theta^{\alpha}\Gamma(\alpha)}e^{-\frac{x}{\theta}} \nu(d\alpha).$$

The restriction of the support of $\nu$ to $\mathbb N$, setting e.g., $\bm \omega = \left\{\omega_i = \nu(i),i \in \mathbb N\right\}$, give rise to the $\texttt{MixedErlang}(\bm \omega,\theta)$ distribution. The number of non-zero $\omega_i$'s is usually finite, and in practice we denote $n = \max_{\omega_i>0} i$. One of the reasons that drove the literature interests upon this class is the denseness results by Tijms~\cite[p. 163]{tijms1994stochastic}. Tijms show that this restricted class can approximate arbitrary distributions on the positive real line. More precisely, if we let $F$ denote a distribution function on the real line, setting weights as $\omega_i = F(i\theta) - F((i-1)\theta)$, we have that the corresponding distribution function converges to $F$ pointwise as $\theta \to 0$, that is 
$$\forall x \in \mathbb R_+, \;F_{\bm \omega,\theta}(x) \xrightarrow[]{\theta \to 0} F(x).$$

A multivariate extension of this result can be found in \cite[Theorem 2.1]{leeModelingDependentRisks2012}. The result can be easily proved using the convergence of moment generating functions. 
The result is constructive, and gives, for any scale $\theta$, a set of weights $\omega_i = F(i\theta) - F((i-1)\theta)$ that, if we could estimate them, would lead to a nice starting point for an iterative algorithm. Note that using them directly as an estimator is not recommended, as, as was already noted by~\cite{leeModelingEvaluatingInsurance2010}, this clearly overfits a given dataset. Moreover, in our case, we do not have enough data do estimate this discretization of $F$. 

Note that, for $X \sim \texttt{MixedGamma}(\nu,\theta)$, $\bm\pi,\bm\mu$ and $\bm\Sigma$ are easily computed. Indeed, it can be checked that binned moments are given by  
\begin{equation}\label{eqn:muboxedmoments}
  \mathbb E\left(X^{k}\mathbb 1_{X \in [a,b[}\right) = \theta^k \int \frac{\Gamma(\alpha+k)}{\Gamma(\alpha)} \left(\gamma\left(\frac{b}{\theta},\alpha+k\right) - \gamma\left(\frac{a}{\theta},\alpha+k\right)\right) \nu(d\alpha),
\end{equation}

where $\gamma(x,\alpha) = \Gamma(\alpha)^{-1}\int_{0}^x t^{\alpha-1} e^{-t} dt$ is the lower incomplete gamma function. The expression in \cref{eqn:muboxedmoments} can be efficiently computed recursively on $k$, and, for  $X\sim \texttt{MixedErlang}(\bm\omega,\theta)$, recursively on atoms of $\nu$. From \cref{eqn:piandmu} we see that these integrals are enough to compute $\bm\pi,\bm\mu$ and $\bm\Sigma$, and from \cref{eqn:llh} our loglikelihood. With these implicit definitions, the loglikelihood of a $\texttt{MixedErlang}(\bm\omega,\theta)$ for our sample is denoted $\ell(\bm\omega,\theta) = \ell(\bm\pi,\bm\mu,\bm\Sigma)$.

Due to overfitting concerns on one side, and due to spiky behavior of Erlang mixtures, we need to include a regularization mechanism in the fitting process. Smoothness in a density estimate is often obtained (see e.g.,\,\cite{guiFittingMultivariateErlang2021a} in the context of Erlang mixtures) by penalizing large norms of the density derivatives through the addition of a penalty $$\widetilde{\mathrm{Pen}_r}(\lambda,f) = \frac{\lambda}{2}\int f^{(r)}(x)^2 dx,$$ for a given positive integer $r$. 

For $f \sim \texttt{MixedErlang}(\bm \omega, \theta)$, $\widetilde{\mathrm{Pen}_r}(\lambda,f) = \frac{\lambda}{2} \theta^{-(2r+1)} \bm \omega'\tilde{\bm P_r}\bm \omega$, where $\tilde{\bm P_r}$ is a fixed positive semidefinite dense matrix.
The elements of the matrix $\tilde{\bm P_r}$ have closed form expressions given by $$\tilde{P}_{r,i,j} = \sum_{k=0}^{r} \sum_{\ell = 0}^{r} c_{r,k} c_{r,\ell} \frac{\Gamma(i+j-k-\ell-1)}{\Gamma(i-k)\Gamma(j-\ell)} 2^{-(i+j-k-\ell-1)} \mathbb 1_{i-k >0, j-l > 0},$$
where $c_{r,k}$ are finite difference coefficients of order $r$, defined recursively by $c_{0,0} = c_{1,0} = 1$ and $c_{r,k} = \left(c_{r-1,k} - c_{r-1,k-1}\right)\mathbb 1_{k \in \left\{0,...r\right\}}$.
While looking for an optimal hyperparameter $\lambda$, the term $\theta^{-(2r-1)}$ could be simply incorporated into $\lambda$ and removed from the equation. However, the matrix $\tilde{\bm P_r}$ is dense, which can be computationally problematic, and in our low information settings we do not have the possibility to run cross-validation strategies to find the right $\lambda$.

We solve these issues by considering instead another penalty, drawn from the P-splines literature. Instead of penalizing the continuous roughness of $f$, we penalize the discrete roughness of the sequence of modes of the underlying Erlang densities. Neglecting a multiplicative $\theta^{-1}$ factor, the mode of the $(i+1)^{\text{th}}$ Erlang density in the mixture is positioned at 
\begin{equation}\label{eq:couples}
  (x_i,y_i) = \left(\theta i, \frac{i^i e^{-i}}{i!}\right).
\end{equation}

Note that $\bm x = \{x_1,x_2,...\}$ are regularly distributed on the positive real line. Therefore, we monitor the regularity of the density estimate by the regularity of the sequence of weighted modes $\omega_1y_1,\omega_2y_2,...$ Let us denote by $\bm D_r$ the $r^{\text{th}}$ difference matrix, of size $(n-r,n)$, defined by 
$$D_{r,k,l} = c_{r,l-k}y_{l}\mathbb 1_{l-k \le r},$$ such that $\bm D_r\omega$ are the finite differences of the sequence of modes. The matrix $\bm D_r'\bm D_r$ is sparse (only its $(2r-1)$ central diagonals are non-zero) which will be computationally efficient. 
The corresponding penalty writes 
$$\frac{\lambda}{2} \lVert \bm D_r\bm \omega \rVert_2^2.$$

We can interpret this penalization as a prior on $r^{\text{th}}$ order differences of coefficients $\bm D_r\bm\omega \vert \lambda \sim \texttt{Normal}(\bm 0, \lambda^{-1} \bm I)$, where $\bm I$ represent the identity matrix. 
Assuming furthermore that we assign a $\texttt{Gamma}(a_{\lambda},b_{\lambda})$ prior on $\lambda$, with hyperparameters $a_{\lambda},b_{\lambda}$ such that the prior is uninformative (high variance), the final penalization is 
$$\mathrm{Pen}_r(\lambda,\bm\omega) = \frac{1}{2}\left\{-(n-r)\log(\lambda) + \lambda\lVert \bm D_r\bm \omega \rVert_2^2 \right\} + \left\{\lambda b_{\lambda}^{-1} - (a_{\lambda}-1)\log(\lambda)\right\}.$$

Thus, we set, still neglecting additive constants, 
\begin{align*}
  \ell(\bm \omega, \theta,\lambda) = &\ell(\bm\omega,\theta) - \mathrm{Pen}_r(\lambda,\bm\omega)\\
  = &\left\{\hat{\bm\pi}'\log(\bm\pi)\right\} - \frac{1}{2} \left\{\log|\bm\Sigma| +\lVert\bm \mu - \hat{\bm \mu}\rVert_{\Sigma}^2\right\} \\
  &- \frac{1}{2}\left\{-(n-r)\log(\lambda) + \lambda\lVert \bm D_r\bm \omega \rVert_2^2 \right\}\\
  &- \left\{\lambda b_{\lambda}^{-1} - (a_{\lambda}-1)\log(\lambda)\right\},
\end{align*}

where the link between the triplet $\left(\bm\pi,\bm\mu,\bm\Sigma\right)$ and parameters $\left(\bm\omega,\theta\right)$ is here implicitly assumed as given by \cref{eqn:piandmu,eqn:muboxedmoments}. 
Note that the added components do not depend on $(\omega,\theta)$ and are thus compatible with other approaches (like cross-validation) that uses the direct penalty $\frac{\lambda}{2} \lVert \bm D_r\bm \omega \rVert_2^2.$ On the other hand, they complete the loglikelihood which allows a Bayesian interpretation.

The marginal loglikelihood for $\lambda$ can be obtained using the following identity: 
\begin{equation}\label{eqn:bayes_ell_lambda}
  \ell(\lambda) = \ell(\bm \omega,\theta,\lambda) - \ell(\bm \omega,\theta | \lambda),
\end{equation}

where $\ell(\bm\omega,\theta | \lambda)$ denotes the conditional log posterior of the Erlang mixtures parameters for a given penalty parameter $\lambda$. 

Consider the following notations for the Hessian matrices involved in our problem, which are linked together by a linearity in $\lambda$: 
\begin{equation*}
  \begin{aligned}
    \bm H(\bm \omega,\theta) &= -\frac{\partial^2}{\partial^2 (\bm\omega,\theta)} \ell(\bm\omega,\theta)\\
    \bm P_r &= \frac{\partial}{\partial\lambda}\frac{\partial^2}{\partial^2 (\bm\omega,\theta)} \mathrm{Pen}_r(\lambda,\bm\omega) = \begin{pmatrix}
      \bm D_r'\bm D_r & 0\\ 
      0 & 0
    \end{pmatrix}\\
    \bm H(\bm\omega,\theta,\lambda) &= -\frac{\partial^2}{\partial^2 (\bm\omega,\theta)} \ell(\bm\omega,\theta,\lambda) = \bm H(\bm \omega,\theta) + \lambda \bm P_r.
  \end{aligned}
\end{equation*}

The conditional log posterior $\ell(\bm \omega,\theta | \lambda)$ can then be approximated using a Laplace approximation. 
The Laplace approximation (see e.g., \cite{tierney1986accurate}) is a computationally efficient method to approximate the posterior distribution of penalized parameters. It has been applied to a wide range of Bayesian models, including generalized linear models \cite{tierney1986accurate}, mixed effects models \cite{breslowApproximateInferenceGeneralized1993}, and latent variable models \cite{rue2005gaussian}.
Several studies have shown that Laplace approximation can provide accurate and efficient estimates for P-spline models, even for large datasets and complex models with high-dimensional spline bases~\cite{fahrmeir2001bayesian,wood2017generalized,lambert2023penalty,gressani2021laplace}. Overall, Laplace approximation remains a popular method for Bayesian inference in the P-splines literature. 

Here, we approximate the posterior distribution around the posterior mode $(\hat{\bm\omega},\hat{\theta})$, that is (neglecting again additive constants):
\begin{equation}\label{eqn:laplace_approx}
  \ell(\bm \omega,\theta | \lambda) \approx \frac{1}{2}\log\lvert \bm H(\hat{\bm\omega},\hat\theta,\lambda) \rvert - \frac{1}{2} \lVert(\bm\omega,\theta) - (\hat{\bm\omega},\hat\theta)\rVert_{\bm H(\hat{\bm\omega},\hat\theta,\lambda)^{-1}}^2.
\end{equation}

This approximation become particularly useful if we note that \cref{eqn:bayes_ell_lambda} is true for any $\bm \omega,\theta$. Thus, by evaluating \cref{eqn:laplace_approx} near the posterior mode $(\hat{\bm\omega},\hat\theta)$, with substitution in \cref{eqn:bayes_ell_lambda}, we obtain: 

$$\ell(\bm \omega,\theta | \lambda) \approx \frac{1}{2}\log\lvert \bm H(\bm\omega,\theta,\lambda) \rvert.$$

Taking the derivative with respect to $\lambda$, using the fact that for a positive-definite matrix $\bm A(x)$, $\frac{\partial}{\partial x} \log \lvert \bm A(x)\rvert = \mathrm{Tr}\left\{A(x)^{-1} \frac{\partial}{\partial x} \bm A(x)\right\}$, one has: 
\begin{equation*}
  \begin{aligned}
    \frac{\partial}{\partial\lambda} \ell(\bm \omega,\theta | \lambda)
    &\approx \frac{\partial}{\partial\lambda} \frac{1}{2}\log\lvert \bm H(\bm\omega,\theta,\lambda) \rvert\\
    &= \frac{1}{2} \mathrm{Tr}\left\{\bm H(\bm\omega,\theta,\lambda)^{-1}\frac{\partial}{\partial\lambda}\bm H(\bm\omega,\theta,\lambda)\right\}\\
    &= \frac{1}{2} \mathrm{Tr}\left\{\bm H(\bm\omega,\theta,\lambda)^{-1}\bm P_r\right\}\\
    &= \frac{1}{2} \mathrm{Tr}\left\{\left(\bm H(\bm\omega,\theta) + \lambda \bm P_r\right)^{-1}\bm P_r\right\},
  \end{aligned}
\end{equation*}

By a recursive application of the Woodbury matrix identity, for two matrices $\bm A$ and $\bm B$ we have that $(\bm A - \bm B)^{-1} = \sum_{k=0}^{\infty} (\bm A^{-1}\bm B)^k \bm A^{-1},$ assuming the sum converges, that is assuming the spectral radius of $\bm A^{-1}\bm B$ is less than one. In our contexts, this yields : 

\begin{equation*}
  \left(\bm H(\bm\omega,\theta) - (-\lambda \bm P_r)\right)^{-1}\bm P_r
  = -\sum_{k=0}^{\infty} \lambda^k \left(-\bm H(\bm\omega,\theta)^{-1} \bm P_r\right)^{k+1}.
\end{equation*}

Assuming that $\eta_1,...\eta_{n+1}$ are eigenvalues of $-\bm H(\bm\omega,\theta)^{-1} \bm P_r$, we therefore have: 

\begin{equation*}
  \begin{aligned}
  \mathrm{Tr}\left\{\left(\bm H(\bm\omega,\theta) + \lambda \bm P_r\right)^{-1}\bm P_r\right\}
  &= -\sum_{k=0}^{\infty} \lambda^k \mathrm{Tr}\left\{\left(-\bm H(\bm\omega,\theta)^{-1} \bm P_r\right)^{k+1}\right\} \\
  &= -\sum_{k=0}^{\infty} \lambda^k \sum_{i=1}^{n+1} \eta_i^{k+1} \\
  &= -\sum_{i=1}^{n+1} \frac{\eta_i}{1 - \lambda\eta_i}.\\
  \end{aligned}
\end{equation*}

Finally,
\begin{align}\label{eqn:lambda_derivative}
    \frac{\partial}{\partial\lambda} \ell(\lambda) 
    &= \frac{\partial}{\partial\lambda} \ell(\bm \omega,\theta,\lambda) - \frac{\partial}{\partial\lambda} \ell(\bm \omega,\theta | \lambda)\nonumber\\
    &= -\frac{1}{2}\left\{ -\frac{n-r}{\lambda} + \lVert \bm D_r'\bm \omega \rVert_2^2\right\} + \left\{\frac{a_{\lambda}-1}{\lambda} - b_{\lambda}^{-1}\right\}  - \frac{1}{2} \mathrm{Tr}\left\{\left(\bm H(\bm\omega,\theta) + \lambda \bm P_r\right)^{-1}\bm P_r\right\}\nonumber\\
    &= \frac{-1}{2}\left\{\lVert \bm D_r'\bm \omega\rVert_2^2 - \frac{(n-r) + 2a_{\lambda}-2}{\lambda} + 2b_{\lambda}^{-1} -\sum_{i=1}^{n+1} \frac{\eta_i}{1 - \lambda\eta_i} \right\},
\end{align}

which, since $\bm\eta$ and $\lVert \bm D_r'\bm \omega \rVert_2^2$ are constants w.r.t. $\lambda$, gives us a simple non-linear (polynomial of degree $n+1$) equation to find critical points of $\ell(\lambda)$, which can be solved with, e.g., a fixed-point method of any gradient ascent method, reaching the closest local maxima.

Thus, \cref{alg:cap} uses an iterative algorithm, by recursively:\begin{enumerate}
  \item Update $(\bm\omega,\theta) = \mathrm{arg}\min \ell(\bm\omega,\theta,\lambda)$ at the current value of $\lambda$,
  \item Update $\lambda = \mathrm{arg}\min \ell(\lambda)$ at the current value of $(\bm\omega,\theta)$ by setting \cref{eqn:lambda_derivative} to $0$,
\end{enumerate}

until certain convergence criterions are met. The iterative procedure would provide, under the Bayesian specification we described, the optimal parameters for our problem.

Assume that the optimization problem for $(\bm\omega,\theta)$ stopped at a mode of the penalized loglikelihood for the selected parameter $\lambda$. The final value of $\bm H(\bm\omega,\theta,\lambda)$ could be used to construct confidence intervals through, again, Laplace approximations. For a real-valued function $(\bm\omega,\bm\theta) \mapsto f(\bm\omega,\bm\theta) \in \mathbb R$, a confidence interval for $f(\bm\omega,\bm\theta)$ can be obtained by a simple Gaussian approximation and the Delta method~\cite{oehlert1992note}.
The Delta method gives an approximation for the variance of $f(\bm\omega,\bm\theta)$ as $N^{-1}\lVert\frac{\partial}{\partial(\bm\omega,\theta)}f(\bm\omega,\theta) \rVert_{\bm H(\bm \omega,\theta,\lambda)}$.
This method provides fast and reliable confidence intervals for any quantity, e.g., for the density itself (pointwise), the quantile function (again, pointwise), or the VaR / TVaR approximations at higher levels.

\section{Numerical examples}\label{sec:numerical_investigation}

In this section, we provide a few numerical examples of our approach. Recall that our final algorithm as the following form :
\begin{algorithm}
  \small
  \DontPrintSemicolon 
  \KwIn{Information about the partition and the moments $\mathcal B, \bm k, \hat{\bm \pi}, \hat{\bm \mu}$.}
  \KwOut{Estimated parameters $\hat{\theta},\hat{\bm\omega}$.}
  Initialize $\lambda = 0$ and $\bm \omega,\theta$ with random (positive) noise.\;
  \While{not converged}{
    $\bm \omega,\theta \gets \arg\min \ell(\bm \omega,\theta,\lambda)$ \tcc*{using e.g., BFGS}
    $\bm H \gets - \frac{\partial^2}{\partial^2 \bm \omega} \ell(\bm \omega, \theta,0)$ \\
    $\eta = \mathrm{eigenvalues}(-\bm H^{-1} \bm P_r)$\\
    $\lambda = \arg\max \ell(\lambda)$ \tcc*{Using \cref{eqn:lambda_derivative}} 
  }
  \Return{$\bm\omega,\theta$}
  \caption{Constrained density estimation using Erlang mixtures.}
  \label{alg:cap}
\end{algorithm}

Note that the computation of the Hessian $\bm H$ can be extracted from the output of the solver used at the previous step, if a second order solver is used. If a first order solver is used, like BFGS, an approximation of the inverse Hessian is usually available at the end of the routine.
In practice, our implementation using the Julia language uses convergence criterions on $\bm\omega,\theta,\lambda$. Moreover, a small $\epsilon$ can be added to the diagonal of $\bm H$ to improve the condition number of the Hessian matrix and facilitate the numerical inversion. Technical implementation relaxes the positivity constraint on the parameters by a square transformation.

\subsection{Lognormal example}

To study the performance and behavior of the proposed scheme, we first revisit a simple example from~\cite{cossette2016moment}, the log-Normal case. Although the log-normal density is unimodal, and our approximation class is clearly targeted at multimodal behaviors, the log-normal density is also heavy tailed and thus a pertinent benchmark for the practitioner. Of course, when the number of atoms of $\nu$ goes to infinity the log-normal is inside our class, and we thus expect good results. A proper analytic expression for $\nu$ can be obtained using the Thorin measure of the log-normal distribution from~\cite{bondessonLevyMeasureLognormal2002,furmanLogNormalConvolutionsAnalyticalNumerical2017,lavernyEstimationMultivariateGeneralized2021,lavernyEstimationHighDimensional2022} and then the mixture representation from~\cite[Theorem 4.1.1]{bondessonGeneralizedGammaConvolutions1992} (or the one from~\cite{moschopoulosDistributionSumIndependent1985}). For practical reasons however, finite mixtures are preferable as they allow computations to be carried out.

\begin{dataset}[Lognormal]\label{ex:lognormal} We sample $N=750$ independent observations from a random variable $X \sim \texttt{LogNormal}(0,0.5)$. The partition is set arbitrarily at the quantiles of the sample at levels $\alpha \in \{0,0.5,0.9,0.99,1\}$, so that it has $J=4$ bins, and we set $\bm k = (4,4,4,1)$ for the number of observed moments within each bin. The corresponding dataset was given in~\cref{tab:exemple_1}. 
\end{dataset}

\begin{figure}[H]
  \centering
  \includegraphics[width=\textwidth]{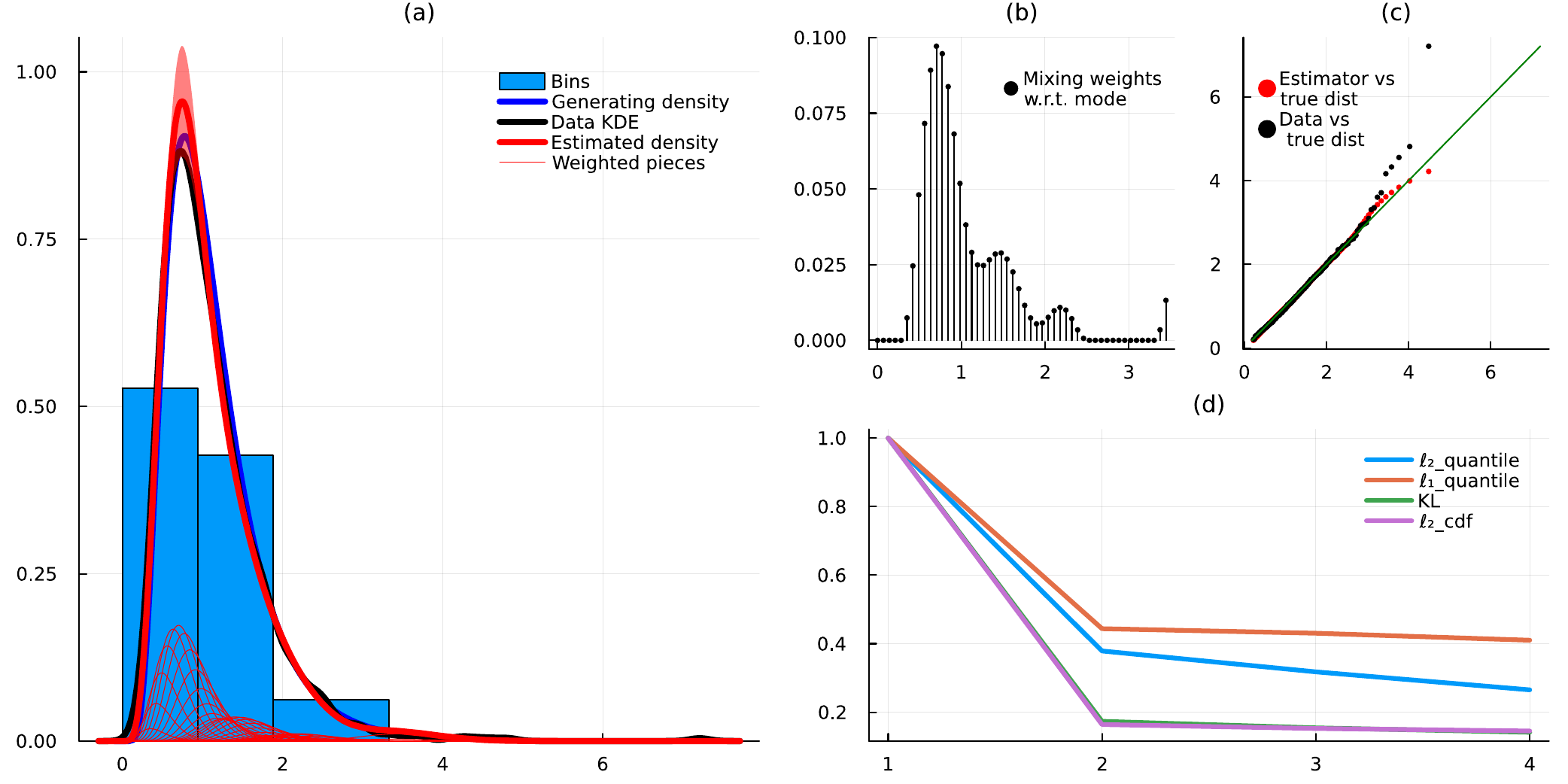}
  \caption{\label{fig:LogNormal:750:statistics_wrt_k_v2} General results on \cref{ex:lognormal}. The histogram in (a) represents $\hat{\boldsymbol \pi}$, while the red curve represents the density estimate and the red region the (pointwise) confidence interval around the estimated density. Couples $(x_i,\omega_iy_i)$ from Equation~\eqref{eq:couples} are represented in panel (b). Panel (c) represents quantile-quantile plots of the estimated density against the underlying raw data and the true density. Panels (d) represents the median of the distance statistics across resamples for an increasing number of observed local moments, renormalized to be $1$ for the least informative setting.}
\end{figure}

The choice of partition in \cref{ex:lognormal} separates the tail from the bulk of the distribution. We run our algorithm with $n=50$ atoms for $\nu$ (that is, $\nu$'s support is contained in $\{1,...,n\}$). The obtained results are pictured in \cref{fig:LogNormal:750:statistics_wrt_k_v2}.

From the set of moments depicted in Table~\ref{tab:exemple_1}, we obtain the optimal Erlang mixture depicted in Figure~\ref{fig:LogNormal:750:statistics_wrt_k_v2} (a, b, c). Figure~\ref{fig:LogNormal:750:statistics_wrt_k_v2} (d) shows the convergence of the procedure w.r.t. the number of empirical moments in the first three bins (1,2,3 or 4) used to estimate the density, as assessed using $S=50$ resamples from the lognormal model. It shows that the distance between the density estimate and the true density (as measured by $\ell_1$ and $\ell_2$ distances between quantile functions and distribution functions, as well as the Kullback-Leibler divergence) decreases significantly with the number of available empirical local moments. The definition of these performance metrics is recalled in \cref{apx:keystats}.

We see on \cref{fig:LogNormal:750:statistics_wrt_k_v2} (a) that the kernel density estimate (black curve) obtained on the raw (non-grouped) data and our estimate based on Erlang mixtures (red curve) from grouped data in \cref{tab:exemple_1} both have troubles reproducing correctly the value at the mode of the true density (blue curve). Our Laplace-approximated confidence point-wise confidence interval confirms the local difficulty. From the quantile-quantile plots on \cref{fig:LogNormal:750:statistics_wrt_k_v2} (c) we see that the fit is satisfactory, except in the far tail (last $5$ to $10$ points over $N=750$). A Kolmogorov-Smirnov test between the true underlying density and the estimated density has a (two-sided) p-value of $\num{0.631}$, failing to reject the null hypothesis that the distributions are the same at that sample size.

To observe the convergence of our estimator with respect to the number of observations, we used $S=50$ resamples from \cref{ex:lognormal} with different number of observations $N$. We report in \cref{tab:LogNormal:quantiles_global} the mean value, bias, standard deviation and RMSE of quantiles (at several levels $\alpha$) obtained from our resamples.

\begin{table}[H]
  \scriptsize
  \caption{Mean, bias, standard deviation and Root Mean Square Error of obtained quantiles at several levels $\alpha$ on resamples from \cref{ex:lognormal}, fitted on several resamples of varying sizes $N$.}
  \label{tab:LogNormal:quantiles_global} 
  \centering
  \begin{tabularx}{\linewidth}{c|XXXXX}
\toprule 
$\alpha$ & 0.5 & 0.9 & 0.95 & 0.99 & 0.995 \\
\midrule 
True values & \num{1.000} & \num{1.898} & \num{2.276} & \num{3.200} & \num{3.625} \\
\midrule 
\multicolumn{1}{c}{Mean} \\
$N = 250$ & \num{1.021} & \num{1.877} & \num{2.307} & \num{3.169} & \num{3.434} \\
$N = 500$ & \num{1.006} & \num{1.909} & \num{2.402} & \num{3.253} & \num{3.497} \\
$N = 750$ & \num{1.000} & \num{1.890} & \num{2.398} & \num{3.276} & \num{3.519} \\
$N = 1000$ & \num{0.997} & \num{1.879} & \num{2.372} & \num{3.252} & \num{3.516} \\
$N = 2000$ & \num{1.002} & \num{1.881} & \num{2.368} & \num{3.247} & \num{3.515} \\
\midrule 
\midrule 
\multicolumn{1}{c}{Std. dev.} \\
$N = 250$ & \num{0.033} & \num{0.112} & \num{0.272} & \num{0.442} & \num{0.469} \\
$N = 500$ & \num{0.021} & \num{0.097} & \num{0.183} & \num{0.298} & \num{0.337} \\
$N = 750$ & \num{0.021} & \num{0.081} & \num{0.185} & \num{0.206} & \num{0.234} \\
$N = 1000$ & \num{0.020} & \num{0.085} & \num{0.193} & \num{0.178} & \num{0.182} \\
$N = 2000$ & \num{0.012} & \num{0.059} & \num{0.179} & \num{0.130} & \num{0.161} \\
\midrule 
\multicolumn{1}{c}{RMSE} \\
$N = 250$ & \num{0.040} & \num{0.114} & \num{0.274} & \num{0.443} & \num{0.507} \\
$N = 500$ & \num{0.022} & \num{0.098} & \num{0.223} & \num{0.303} & \num{0.361} \\
$N = 750$ & \num{0.021} & \num{0.081} & \num{0.222} & \num{0.219} & \num{0.257} \\
$N = 1000$ & \num{0.020} & \num{0.087} & \num{0.215} & \num{0.186} & \num{0.212} \\
$N = 2000$ & \num{0.012} & \num{0.061} & \num{0.201} & \num{0.138} & \num{0.195} \\
\bottomrule 
\end{tabularx}

\end{table}

We note on \cref{tab:LogNormal:quantiles_global}, especially on the root-mean-square errors, that all quantiles seem convergent with the number of observations. In particular, quantiles higher than the last bin edge are also more and more precisely estimated with the increasing number of observations. 

In \cref{fig:LogNormal:750:statistics_wrt_k_v2} (d), a few key statistics (detailed in \cref{apx:keystats}) are averaged on the resamples for several values of $\bm k$ (but fixed number of observations $N$). This plot shows that increasing the number of moments in each bin does have an impact on the quality of the estimator, at least on this example. However, the limit is not exactly the lognormal distribution, as the log-normal distributions cannot be characterized by its moments.

\begin{figure}[H]
  \centering
  \begin{subfigure}{0.45\textwidth}
    \includegraphics[width=\textwidth]{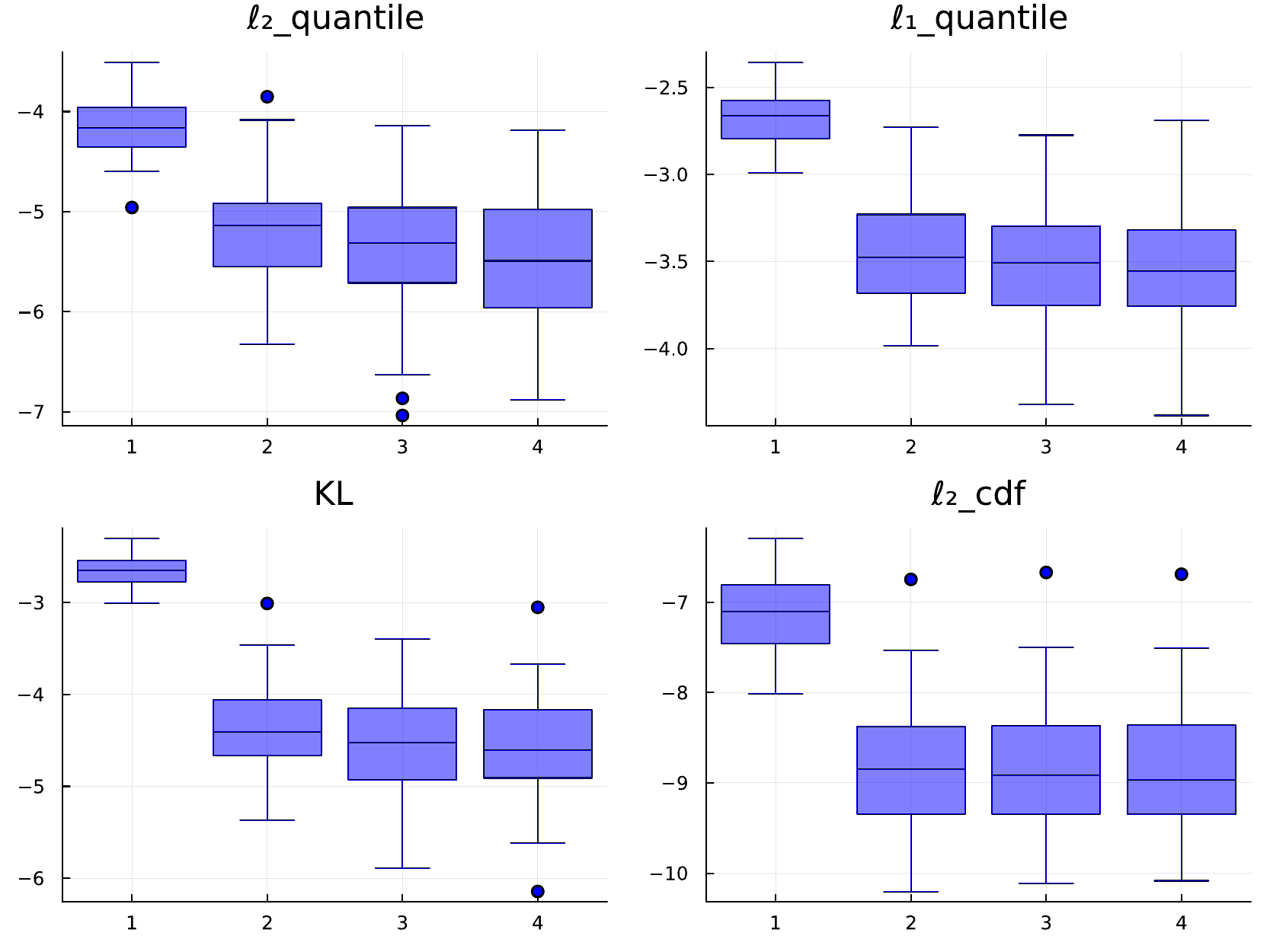}
    \caption{\label{fig:LogNormal:750:boxplots_log}}
  \end{subfigure}
  \hfill
  \begin{subfigure}{0.45\textwidth}
    \includegraphics[width=\textwidth]{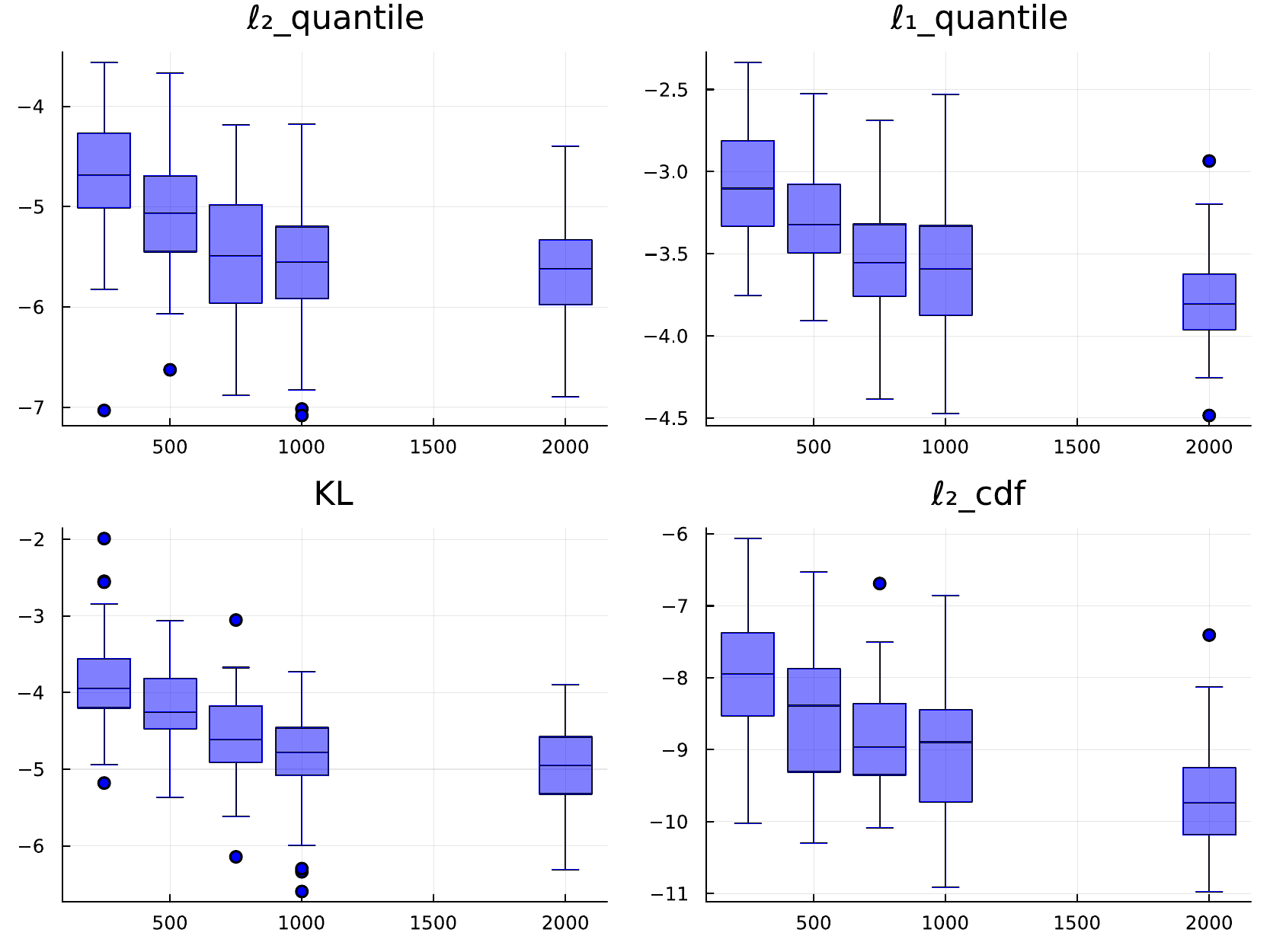}
    \caption{\label{fig:LogNormal:compare_boxplots}}
  \end{subfigure}
  \caption{Comparison of the bootstrap statistics obtained on \cref{ex:lognormal}. Panel (a) : for an increasing number of moments with $\bm k \in \{(1,1,1,1),(2,2,2,1),(3,3,3,1),(4,4,4,1)\}$ and a fixed number of observations $N$. Panel (b): for an increasing number of observations $N \in \{250,500,750,1000,2000\}$ and fixed $\bm k=(4,4,4,1)$.}
\end{figure}

Box plots of these key statistics, for a given value of $\bm k$ but several values of $N$ across the resamples are available in \cref{fig:LogNormal:compare_boxplots}, with clear indications of convergence.

While this example is overall conclusive on the quality of the estimator, we propose a few more examples with additional specific features (such as multimodality) for the target density.

\subsection{A mixture of gamma distributions}

\begin{dataset}[Mixture of gammas]\label{ex:MixGamma1} We sample $N=2000$ independent observations from a random variable $X$ that follows a $\frac{1}{3}$ — $\frac{2}{3}$ mixture of $\texttt{Gamma}(30,1)$ and $\texttt{Gamma}(7,1)$. The partition is set at the quantiles of the sample at levels $\alpha \in \{0,0.5,0.9,0.99,1\}$, so that it has $J=4$ bins, and we use $\bm k = (4,4,4,1)$ for the number of observed moments within each bin.
\end{dataset}

\begin{figure}[H]
  \centering
  \includegraphics[width=\textwidth]{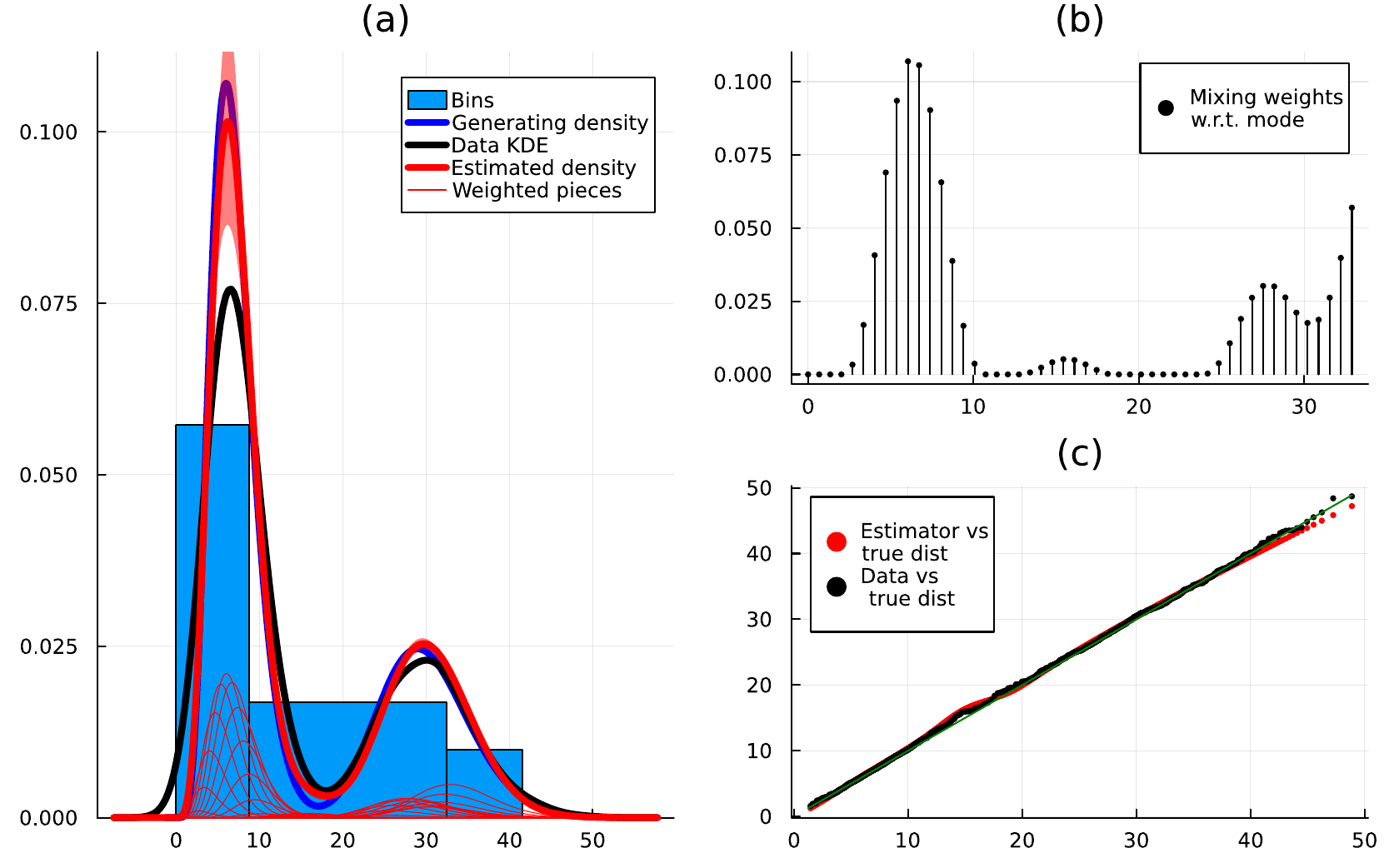}
  \caption{\label{fig:MixGamma1:first_plot} General results on \cref{ex:MixGamma1}. The histogram in (a) represents $\hat{\boldsymbol \pi}$, while the red curve represents the density estimate and the red region the (pointwise) confidence interval around the estimated density. Couples $(x_i,\omega_iy_i)$ from Equation~\eqref{eq:couples} are represented in panel (b). Panel (c) represents quantile-quantile plots of the estimated density against the underlying raw data and the true density.}
\end{figure}

Our estimator should be able to recover the features of the underlying distribution fairly correctly since the true distribution actually belongs to the class of Erlang mixtures. However, identifying the two atoms of the underlying measure is not possible in our $L_2$-penalized settings\footnote{Finding a few atoms would require instead some kind of $L_1$ penalization, see e.g., \cite{yinEfficientEstimationErlang2016,yinConsistencyPenalizedMLEs2019} and their iSCAD penalty}, as confirmed by the results in \cref{fig:MixGamma1:first_plot}. 

We observe a particular difficulty in the estimation of the first mode of the density on \cref{fig:MixGamma1:first_plot}, but our estimator performs better than the standard Gaussian kernel — that has all the observation available and not only the boxed moments. On the other hand, Kolmogorov-Smirnov tests are conclusive (p-value of $\num{0.157}$) on this example.

\subsection{A more challenging example}

To confront our estimator to a dataset that is obviously outside its range, we construct \cref{ex:GaussRevGamma500}.

\begin{dataset}[Mixture of a Gaussian and a (reversed) Gamma~\cite{lambertNonparametricDensityEstimation2023}]\label{ex:GaussRevGamma500} We sample $N=500$ independent observations from a random variable $X$ that follows a 20/80 mixture of $\texttt{Normal}(1,1/3)$ and a reversed $\texttt{Gamma}(11,1/6)$. The partition is set arbitrarily at the quantiles of the sample at levels $\alpha \in \{0,0.5,0.9,0.99,1\}$, so that it has $J=4$ bins, and we use $\bm k = (4,4,4,1)$ for the number of observed moments in each bin.
\end{dataset}

\begin{figure}[H]
  \centering
  \includegraphics[width=\textwidth]{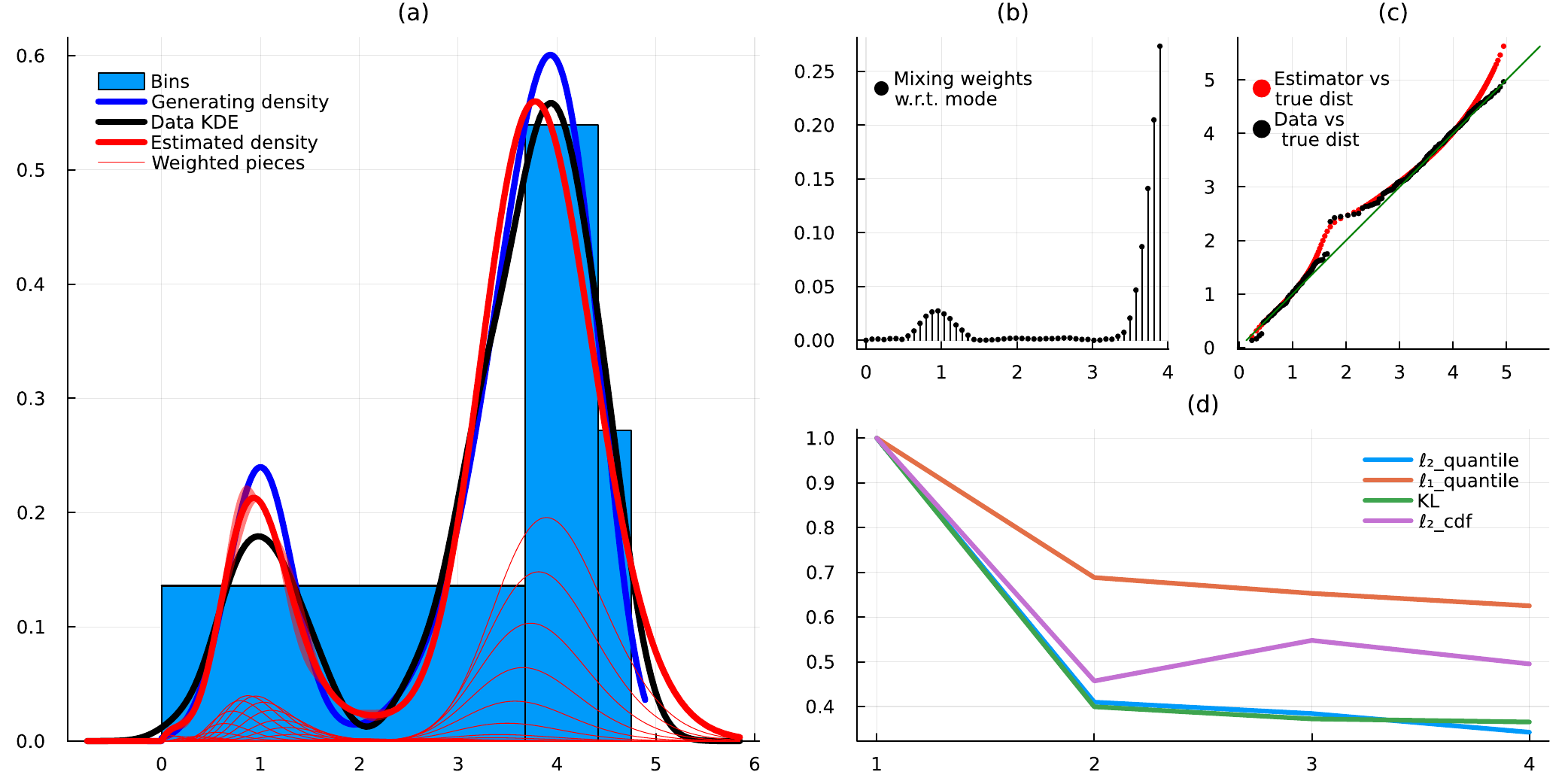}
  \caption{\label{fig:GaussRevGamma:500:statistics_wrt_k_v2} General results on \cref{ex:GaussRevGamma500}. The histogram in (a) represents $\hat{\boldsymbol \pi}$, while the red curve represents the density estimate and the red region the (pointwise) confidence interval around the estimated density. Couples $(x_i,\omega_iy_i)$ from Equation~\eqref{eq:couples} are represented in panel (b). Panel (c) represents quantile-quantile plots of the estimated density against the underlying raw data and the true density. Panels (d) represents the median of the distance statistics across resamples for an increasing number of observed local moments, renormalized to be $1$ for the least informative setting.}
\end{figure}

\begin{figure}[H]
  \centering
  \begin{subfigure}{0.45\textwidth}
    \includegraphics[width=\textwidth]{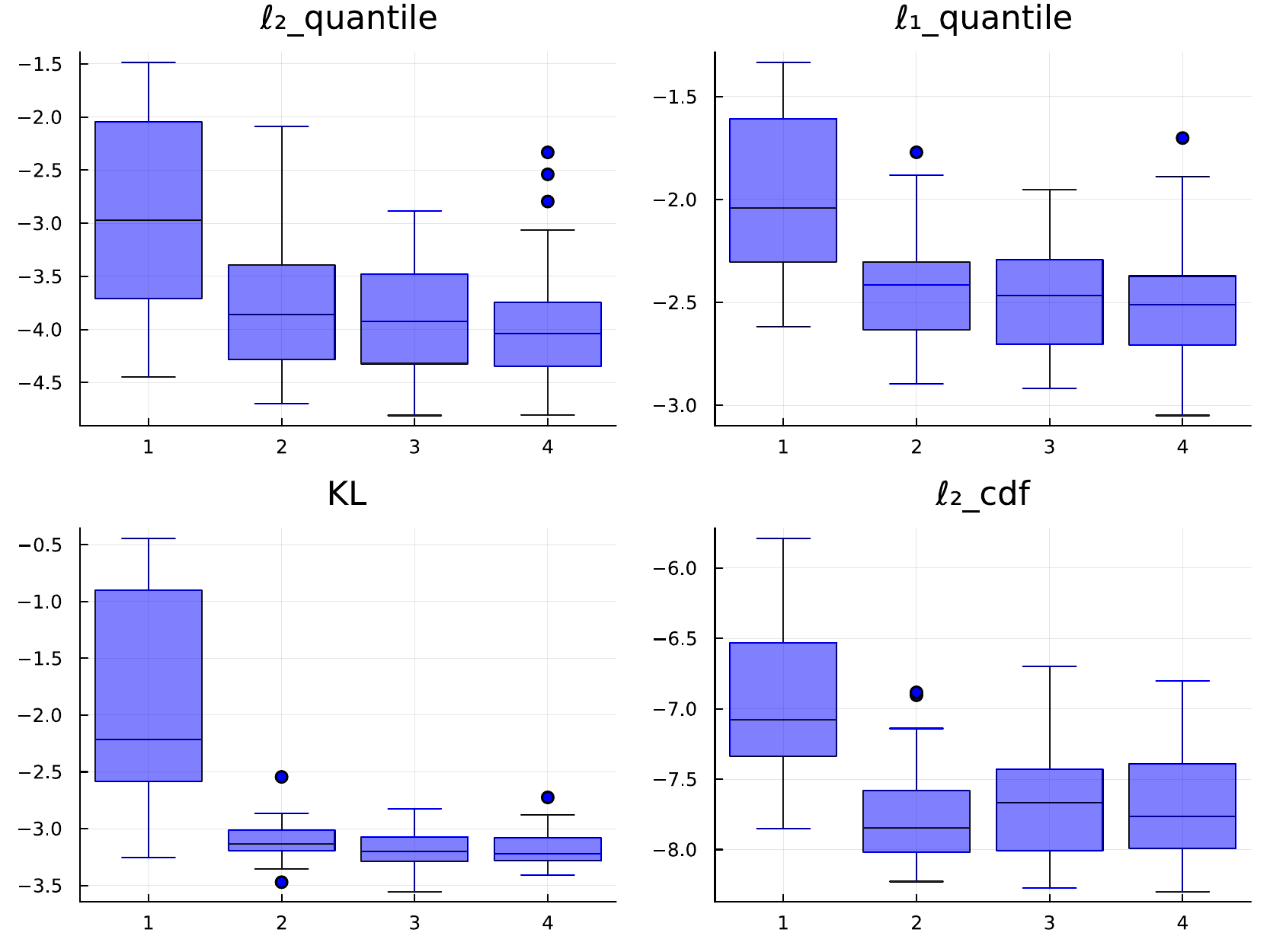}
     \caption{\label{fig:GaussRevGamma:500:boxplots_log}}
  \end{subfigure}
  \hfill
  \begin{subfigure}{0.45\textwidth}
    \includegraphics[width=\textwidth]{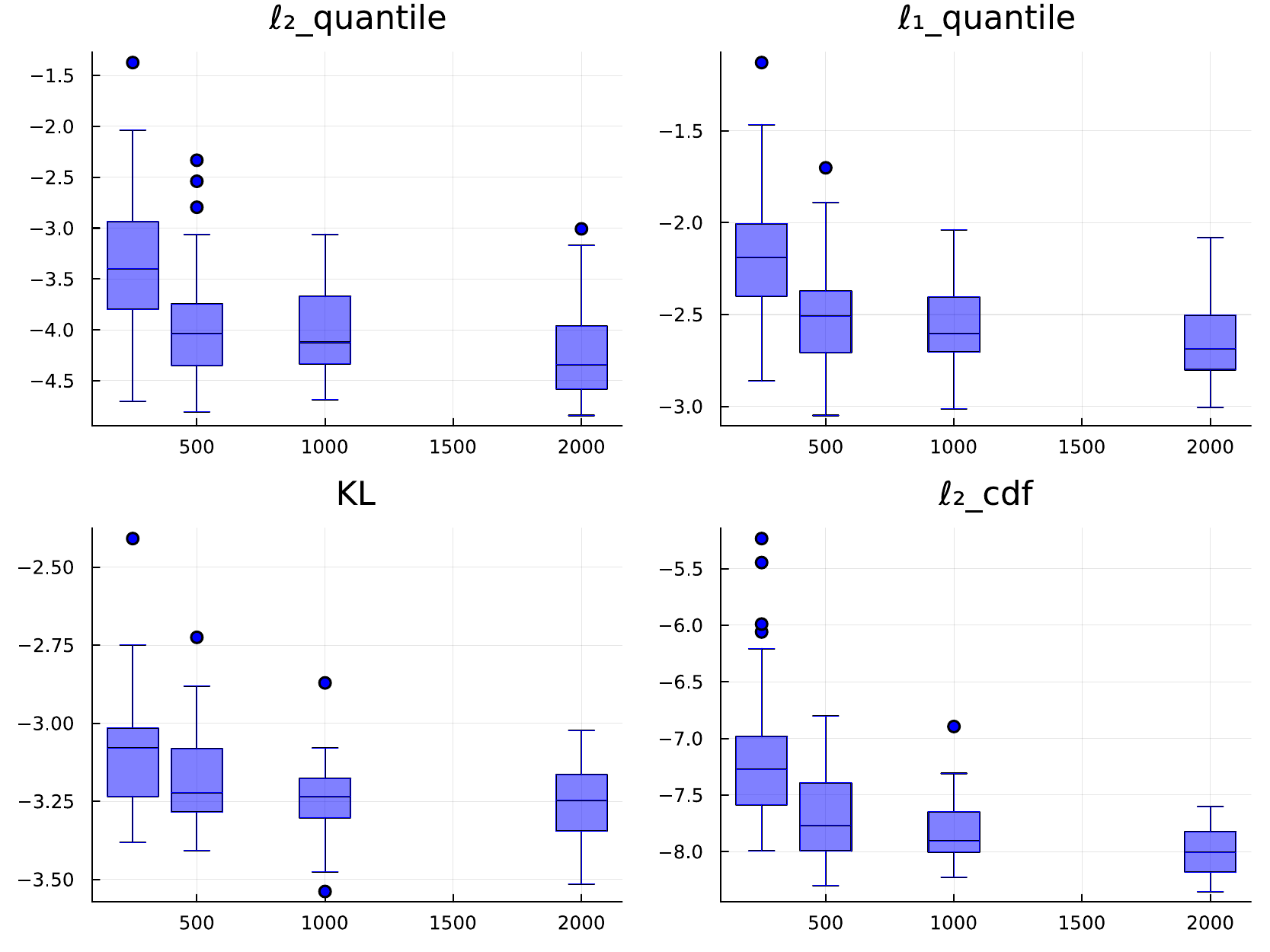}
    \caption{\label{fig:GaussRevGamma:compare_boxplots}}
  \end{subfigure}
  \caption{Comparison of the bootstrap statistics obtained on \cref{ex:GaussRevGamma500}. Panel (a) : for several parameters $\bm k \in \{(1,1,1,1),(2,2,2,1),(3,3,3,1),(4,4,4,1)\}$ and fixed number of observations $N$. Panel (b): for several number of observations $N \in \{250,500,1000,2000\}$ and fixed $\bm k=(4,4,4,1)$.}
\end{figure}

\begin{table}[H]
  \scriptsize
  \caption{Mean, bias, standard deviation and Root Mean Square Error of obtained quantiles at several levels $\alpha$ on resamples from \cref{ex:lognormal}, fitted on several resamples of varying sizes $N$ on \cref{ex:GaussRevGamma500}.}
  \label{tab:GaussRevGamma:quantiles_global} 
  \centering
  \begin{tabularx}{\linewidth}{c|XXXXX}
\toprule 
$\alpha$ & 0.5 & 0.9 & 0.95 & 0.99 & 0.995 \\
\midrule 
True values & \num{3.643} & \num{4.376} & \num{4.530} & \num{4.778} & \num{4.857} \\
\midrule 
\multicolumn{1}{c}{Mean} \\
$N = 250$ & \num{3.589} & \num{4.469} & \num{4.715} & \num{5.189} & \num{5.366} \\
$N = 500$ & \num{3.606} & \num{4.483} & \num{4.729} & \num{5.202} & \num{5.378} \\
$N = 1000$ & \num{3.612} & \num{4.490} & \num{4.737} & \num{5.211} & \num{5.390} \\
$N = 2000$ & \num{3.610} & \num{4.480} & \num{4.726} & \num{5.199} & \num{5.377} \\
\midrule 
\midrule 
\multicolumn{1}{c}{Std. dev.} \\
$N = 250$ & \num{0.054} & \num{0.048} & \num{0.050} & \num{0.056} & \num{0.059} \\
$N = 500$ & \num{0.039} & \num{0.042} & \num{0.045} & \num{0.054} & \num{0.058} \\
$N = 1000$ & \num{0.028} & \num{0.032} & \num{0.036} & \num{0.044} & \num{0.045} \\
$N = 2000$ & \num{0.014} & \num{0.023} & \num{0.027} & \num{0.034} & \num{0.038} \\
\midrule 
\multicolumn{1}{c}{RMSE} \\
$N = 250$ & \num{0.077} & \num{0.105} & \num{0.191} & \num{0.414} & \num{0.513} \\
$N = 500$ & \num{0.054} & \num{0.115} & \num{0.204} & \num{0.427} & \num{0.525} \\
$N = 1000$ & \num{0.042} & \num{0.119} & \num{0.210} & \num{0.435} & \num{0.535} \\
$N = 2000$ & \num{0.037} & \num{0.107} & \num{0.198} & \num{0.422} & \num{0.522} \\
\bottomrule 
\end{tabularx}

\end{table}

We observe on \cref{tab:GaussRevGamma:quantiles_global} that the bootstrap standard error of the quantiles' estimation is decreasing with the number of observations, but that the root-mean-square error for extreme quantiles of the target distribution, which is based on a comparison of the estimated quantity with its true value, does not. Hence, there is here an error that does not decrease with the number of observations for the estimation of extreme quantiles, see on \cref{fig:GaussRevGamma:500:statistics_wrt_k_v2}, panel (a). Still, generic distances from \cref{fig:GaussRevGamma:compare_boxplots} are nevertheless correctly decreasing. The Kolmogorov-Smirnov p-value of our estimator on \cref{ex:GaussRevGamma500} with $N=500$ is at $\num{0.048}$.

\subsection{A mixture of a (truncated) normal and a beta random variables.}

\begin{dataset}[Mixture of Gaussian and Beta]\label{ex:MixNormalBeta} We sample $N=750$ independent observations from a random variable $X$ that follows a 20/80 mixture of $\texttt{Normal}(1,1/3)$ and $\texttt{Beta}(4,8)$. The partition is set arbitrarily as the quantiles of the sample at levels $\alpha \in \{0,0.5,0.9,0.99,1\}$, so that it has $J=4$ bins, and we use $\bm k = (4,4,4,1)$ for the number of observed moments within each bin.
\end{dataset}

\begin{figure}[H]
  \centering
  \includegraphics[width=\textwidth]{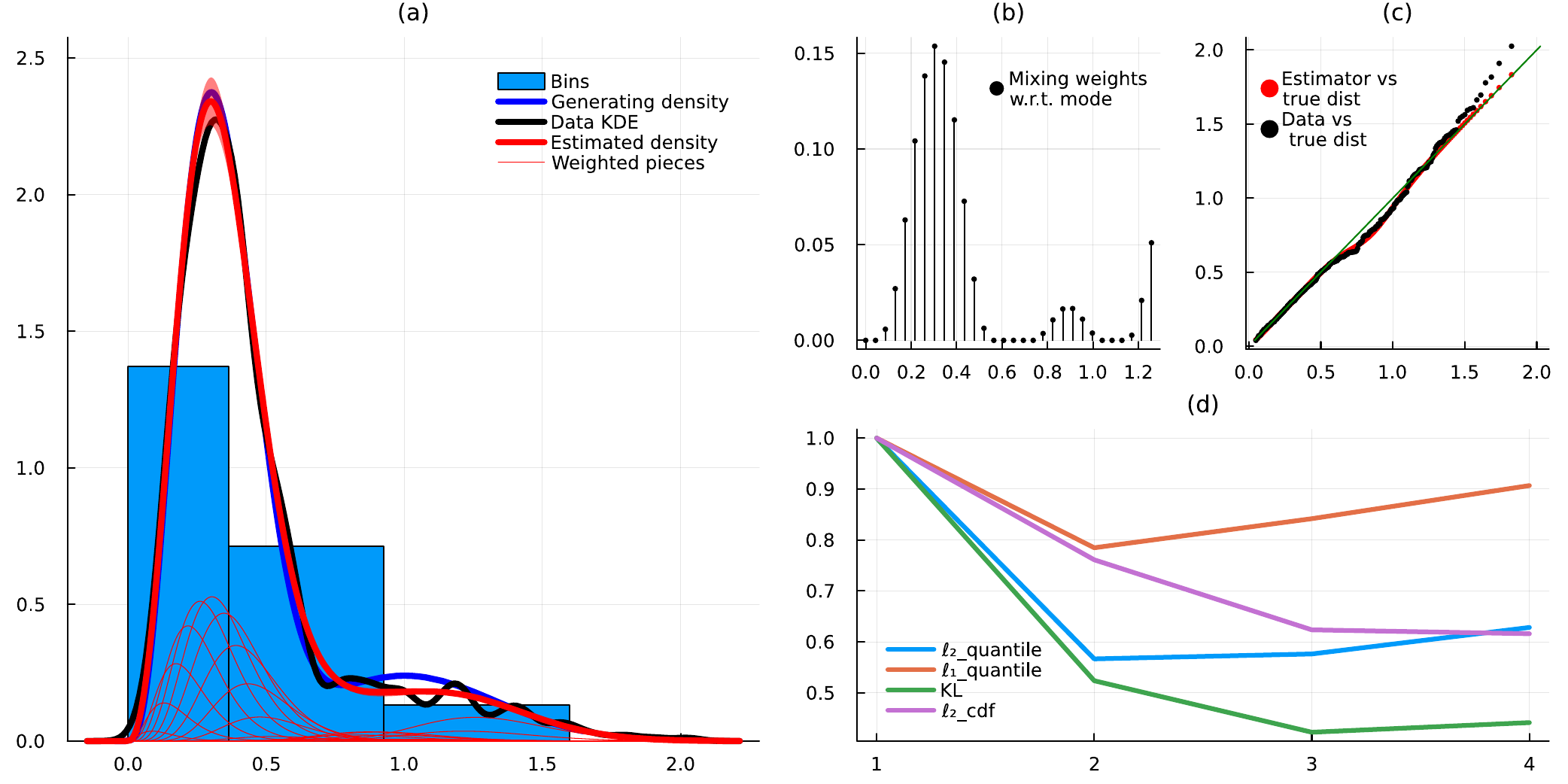}
  \caption{\label{fig:MixNormalBeta:statistics_wrt_k_v2} General results on \cref{ex:MixNormalBeta}. The histogram in (a) represents $\hat{\boldsymbol \pi}$, while the red curve represents the density estimate and the red region the (pointwise) confidence interval around the estimated density. Couples $(x_i,\omega_iy_i)$ from Equation~\eqref{eq:couples} are represented in panel (b). Panel (c) represents quantile-quantile plots of the estimated density against the underlying raw data and the true density. Panels (d) represents the median of the distance statistics across resamples for an increasing number of observed local moments, renormalized to be $1$ for the least informative setting.}
\end{figure}

\Cref{ex:MixNormalBeta}, and in particular panel (a) of \cref{fig:MixNormalBeta:statistics_wrt_k_v2} highlights the smoothing effect of the $L_2$ penalization. Indeed, by comparing the black curve (the Gaussian kernel that has access to (non-grouped) raw data) with the red curve (our estimator), we clearly see the impact of smoothing. The statistics from panel (d) are not obvious to interpret, but the p-value of the Kolmogorov-Smirnov test of $\num{0.306}$ confirms the quality of the estimated density.

\section{Discussion}\label{sec:conclusion}

Boxed moments problems on the form of \cref{tab:exemple_1} may arise for different reasons. Confidentiality issues, for example, is particularly present in the insurance and reinsurance sector. The boxed moments structure is quite permissive as it, e.g., allows to provide a few moments for the \emph{attritional} part, the Value-at-Risk, and potentially the Tail-Value-at-Risk, which all can generally be obtained directly for already-present databases. The boxing also improves the quality of the transmitted information without increasing the moment's order, which could lead to numerical overflows and/or underflows (see e.g., \cite{cossette2016moment}).

If \cite{lambertNonparametricDensityEstimation2023} dealt with the bounded case, we here provide a new estimator that allows for estimation on the (positive) real line. The class of Erlang mixtures was a natural pick due to their denseness into the class of distributions supported on $\mathbb R_+$, and their previous use is insurance losses modeling. These distributions provide closed-form boxed moments, a characteristic that contributes to the performance of our Julian implementation. 

The very nature of our estimation problem requires regularization. In a setting where cross-validation is not possible, we adapted Bayesian techniques from P-splines discrete penalization. We derived an explicit algorithm that leverages this Bayesian nature to automatically tune the penalty, leaving less hyperparameters to the user. Of course, the number of Erlang distributions in the mixtures (our parameter $n$) is another important hyperparameter, but due to the denseness of the class the simple rule of \emph{bigger $n$ is better} can be used.

Our approach could be easily generalized to multivariate settings for estimation on $\mathbb R_+^d$. For that, we could leverage the approach from \cite{yinEfficientEstimationErlang2016,yinConsistencyPenalizedMLEs2019}, granted in quite different settings ($L_1$-penalized settings on full datasets instead of $L_2$-penalized settings on restricted information). Moreover, we think that concurrent $L_1$ and $L_2$ penalization could be beneficial to the sparsity \emph{and} the smoothness of the obtained estimator, if both criterions are of interests to the user. 

\section{Acknowledgements}

Authors acknowledge the support of the ARC project IMAL (grant 20/25-107) financed by the Wallonia-Brussels Federation and granted by the Académie Universitaire Louvain.

\appendix

\section{Commonly used probability distributions}\label{apx:keydists}

We define here the probability distributions used in this paper and there parametrization by providing a characterizing function for each of them: 

\begin{itemize}
  \item $\texttt{Normal}(\bm \mu,\bm \Sigma)$ has density $f(x) = \frac{1}{\sqrt{2\pi\sigma^2}} e^{-\frac{1}{2}\left(\frac{x - \mu}{\sigma}\right)^2}$, for all $\mu \in \mathbb R, \sigma \in \mathbb R_+$.
  \item $\texttt{LogNormal}(\mu,\sigma)$ has density $f(x) = \frac{1}{x \sqrt{2 \pi \sigma^2}} \exp \left( - \frac{(\log(x) - \mu)^2}{2 \sigma^2} \right),  \quad x > 0,$ for all $\mu \in \mathbb R, \sigma \in \mathbb R_+$.
  \item $\texttt{Gamma}(\alpha,s)$ has density $f(x) = \frac{1}{\Gamma(\alpha) s^\alpha }x^{\alpha-1}e^{-\frac{x}{s}}$, for all $\alpha,s \in  \mathbb R_+$ and $\Gamma$ is the gamma function.
  \item $\texttt{Beta}(a,b)$ has density $f(x) = \frac{x^{a-1} (1-x)^{b-1}}{B(a,b)}$, where $a,b \in \mathbb R_+$ and $B$ is the Beta function.
  \item $\texttt{Multinomial}(\bm \pi, N)$ has probability mass function $p(\bm n) = \frac{n!}{n_1!,...,n_m!} \bm\pi^{\bm n}$, where $\bm\pi$ is a probability vector and $\bm n$ is a vector of integers summing to $N$. 
  \item $\texttt{MixedGamma}(\nu, s)$ has density $f(x) = \int f_{\texttt{Gamma}(\alpha,s)}(x) \nu(d\alpha)$ for $\nu \in \mathcal M_+(\mathbb R_+)$ and $s \in \mathbb R_+$.
  \item $\texttt{MixedErlang}(\bm \omega, s) = \texttt{MixedGamma}(\sum_{k=1}^n \omega_k \delta_k, s)$, where $n$ is the length of the probability vector $\omega$.
\end{itemize}

Loglikelihood of these distributions is derived simply by letting $\ell = \log \circ f$. More details on these distributions can be found in the literature.

\section{Key evaluation statistics}\label{apx:keystats}

The paper leveraged a few key statistics to visualize the quality of the resampled models:
\begin{itemize}
  \item The $\ell_2$ squared distance between quantile function $q_1$ and $q_2$: $$\ell_{2}(q_1,q_2) = \int_{0}^1 (q_1(x) - q_2(x))^2 dx.$$
  \item The $\ell_2$ squared distance between distribution function $F_1$ and $F_2$: $$\ell_{2}(F_1,F_2) = \int (F_1(x) - F_2(x))^2 dx.$$
  \item The $\ell_1$ distance between quantile function $q_1$ and $q_2$ (it is technically equal to the $\ell_1$ distance between distribution functions $F_1$ and $F_2$.): $$\ell_{1}(q_1,q_2) = \int_{0}^1 \mathrm{abs}(q_1(x) - q_2(x)) dx.$$ 
  \item The Kullback-Leibler divergence between density functions $f_1,f_2$: $$KL(f_1,f_2) = \int f_1(x)\ln\left(\frac{f_1(x)}{f_2(x)}\right) dx.$$
\end{itemize}

\printbibliography
\end{document}